\documentclass[12pt,reqno]{amsart}
\usepackage{amssymb}
\usepackage{amsfonts}
\usepackage{color}
\usepackage{hyperref}
\usepackage{amsaddr}

\newtheorem{theorem}{Theorem}[section]

\def\F{\mathbb{F}}
\def\T{\mathcal{T}}
\def\K{\mathcal{K}}
\def\D{\mathcal{D}}
\def\N{\mathcal{N}}
\DeclareMathOperator{\Aut}{Aut}

\begin{document}

\title{Some new Steiner designs $S(2,6,91)$}

\author[Kiermaier, Kr\v{c}adinac, Tonchev, Vlahovi\'{c} Kruc, Wassermann]{Michael Kiermaier$^1$,
Vedran Kr\v{c}adinac$^2$,\\Vladimir D.\ Tonchev$^3$, Renata Vlahovi\'{c} Kruc$^4$\\and Alfred Wassermann$^1$}

\address{$^1$Department of Mathematics, University of Bayreuth,\\ 95440 Bayreuth, Germany}

\address{$^2$Faculty of Science, University of Zagreb,\\ Bijeni\v{c}ka cesta~30, 10000 Zagreb, Croatia}

\address{$^3$Mathematical Sciences, Michigan Technological University,\\ Houghton, Michigan, USA}

\address{$^4$Faculty of Teacher Education, University of Zagreb,\\ Savska cesta 77, 10000 Zagreb, Croatia}

\email[Michael Kiermaier]{michael.kiermaier@uni-bayreuth.de}
\email[Vedran Kr\v{c}adinac]{vedran.krcadinac@math.hr}
\email[Vladimir D.\ Tonchev]{tonchev@mtu.edu}
\email[Renata Vlahovi\'{c} Kruc]{renata.vlahovic.kruc@ufzg.hr}
\email[Alfred Wassermann]{alfred.wassermann@uni-bayreuth.de}



\date{May 17, 2025}

\begin{abstract}
The Kramer-Mesner method for constructing designs with a
prescribed automorphism group~$G$ has proven effective
many times. In the special case of Steiner designs, the task
reduces to solving an exact cover problem, with the advantage
that fast backtracking solvers like Donald Knuth's \emph{dancing
links} and \emph{dancing cells} can be used. We find ways to
encode the inherent symmetry of the problem space, induced by
the action of the normalizer of $G$, into a single instance of
the exact cover problem. This eliminates redundant computations
of certain isomorphic search branches, while preventing the
overhead caused by repeatedly restarting the solver.

Our improved approach is applied to the parameters $S(2,6,91)$.
Previously, only four such Steiner designs were known, all of
which had been constructed as cyclic designs over four decades
ago. We find $23$ new designs, each with full automorphism
group of order~$84$.
\end{abstract}

\maketitle

\section{Introduction}\label{intro}

Let $V=\{1,\ldots,v\}$ be a set of $v$ \emph{points}.
A \emph{Steiner design} $S(t,k,v)$ is a set~$\D$ of $k$-subsets
of~$V$, called \emph{blocks}, such that every $t$-subset of~$V$
is contained in exactly one block. A more general concept are
$t$-$(v,k,\lambda)$ \emph{designs}, where every $t$-subset of
points is contained in exactly~$\lambda$ blocks.

We will mainly consider the case $t=2$ and $\lambda=1$, i.e.\
Steiner $2$-designs. In this case, the necessary existence
conditions -- namely the divisibility criterion and Fisher's
inequality -- can be written as
\[
	k-1 \mid v-1\text{,}\quad
	k(k-1) \mid v(v-1)\quad\text{and}\quad
	v\ge k^2-k+1\text{.}
\]
For the general form of these conditions and for other results
and definitions about combinatorial designs, see \cite{BJL99,
CD07, VDT88}.

If we set $k=3$, the conditions become $v\equiv 1, 3 \pmod{6}$
and $v\ge 7$. Kirkman~\cite{TPK1847} famously proved that this is
already sufficient for the existence of Steiner triple systems,
even before Steiner posed the question in~\cite{JS1853}.
Notably, Steiner's general question was not about designs,
although special cases are equivalent to $S(2,3,v)$ and
$S(3,4,v)$.  For an interesting historical account,
see~\cite{CKO23}.  We use the traditional terminology and
attribute all $t$-designs with $\lambda=1$ to Steiner, as defined
in the first paragraph. Regarding $S(2,k,v)$ designs with $k=4$,
the necessary existence conditions are $v\equiv 1, 4 \pmod{12}$,
$v\ge 13$, and for $k=5$ they are $v\equiv 1, 5 \pmod{20}$,
$v\ge 21$. Hanani~\cite{HH75} proved that these conditions are
also sufficient.

For $k=6$, the necessary conditions
$v\equiv 1, 6 \pmod{15}$, $v\ge 31$ are not sufficient.
It is well known that an affine plane of order~$6$,
i.e.\ a $S(2,6,36)$ design, does not exist. Furthermore,
in~\cite{HTJL01} it was proved that there are no
$S(2,6,46)$ designs. According to \cite[Table 3.4]{AG07},
there are $29$ values of~$v$ for which the existence of
$S(2,6,v)$ designs is an open problem, starting with $v=51$,
$61$, $81$, $166\ldots$ and ending with $v=801$. As far as
we know, all these values are still open.

The parameters $S(2,6,91)$ have attracted considerable attention.
The first example of a design with these parameters was found by
Mills in 1975~\cite{WHM75}. A few years later, C.~J.~Colbourn and
M.~J.~Colbourn constructed three further examples~\cite{CC80,
CC82}. Until now, these were the only known $S(2,6,91)$ designs.
The full automorphism groups of the four designs were
calculated in~\cite{ST88}. The design with the largest
automorphism group was named after Gordon McCalla in~\cite{CC82}.
The order of its automorphism group attains an upper bound proved
by Camina and Di Martino~\cite{CDM89}.  Janko and
Tonchev~\cite{JT91} proved that a cyclic $S(2,6,91)$ with an
automorphism group of order greater than $91$ must be one of the
known designs. Baicheva and Topalova~\cite{BT12, BT19} proved
that these four designs are the only cyclic $S(2,6,91)$s.
Recently, Crnković and Dumičić Danilović~\cite{CDD25} classified
all $S(2,6,91)$s with a non-abelian group of order $26$ and proved
that, in this case, two of the known designs are the only examples.

All four known $S(2,6,91)$ designs are cyclic, i.e.\ have an
automorphism of order $91$. Base blocks are given
in~\cite{CDD25}, and in Table~\ref{tab1} we give information
about their full automorphism groups. The designs are named in
accordance with~\cite{CDD25}.

\begin{table}[ht]
\caption{The four cyclic $S(2,6,91)$ designs.}
\label{tab1}
\begin{tabular}{cccc}
\hline
Label & Name & $\Aut$ & $|\Aut|$ \\
\hline
$\D_1$ & Mills & $C_{91} \rtimes C_3$ & $273$ \\
$\D_2$ & McCalla & $C_{91} \rtimes C_{12}$ & $1092$ \\
$\D_3$ & Colbourn & $C_{91} \rtimes C_4$ & $364$ \\
$\D_4$ & Colbourn & $C_{91}$ & $91$ \\
\end{tabular}
\end{table}

In~\cite{CDD25}, the question was raised whether there exist
$S(2,6,91)$ designs with full automorphism groups of order
strictly smaller than $91$. In this paper, we provide an
affirmative answer by proving the following theorem.

\begin{theorem}\label{tm1}
Up to isomorphism, there are exactly $24$ Steiner designs $S(2,6,91)$ with a
group of automorphisms of order $84$ splitting the point set into an
orbit of size~$7$ and an orbit of size~$84$. One of these designs is the
McCalla design $\D_2$. The remaining $23$ designs have no other automorphisms,
i.e.\ their full automorphism groups are of order~$84$.
\end{theorem}

The layout of our paper is as follows. In Section~\ref{sec2}, we describe
general computational methods for the construction of designs with prescribed
groups of automorphisms. In contrast to~\cite{CDD25}, where the method of
tactical decompositions was used, we refine the Kramer-Mesner method for the
construction of Steiner designs. As an illustration, we classify cyclic
$S(2,6,91)$ designs. In Section~\ref{sec3}, we describe the calculations
proving Theorem~\ref{tm1}. Finally, in Section~\ref{sec4}, we give some
concluding remarks and discuss possible directions for further research.

While our goal is to adapt the construction method to the case of Steiner
designs, we aim to keep the algorithms and their implementations as general as
possible. In doing so, we obtain highly reusable algorithms applicable to a broad
class of problems, which also provides a large set of test cases for verifying
the program’s correctness. That said, it should not go unmentioned that, for
fixed parameters and prescribed groups, it may be possible to create customized
algorithms that improve the running times reported in this article -- for example,
by using optimized data structures tailored to the specific situation.

\section{The Kramer-Mesner method for Steiner designs}\label{sec2}

A group~$G$ of automorphisms of a design~$\D$ partitions the sets of
points and blocks into orbits. These orbits form a tactical decomposition
and give rise to the so-called \emph{orbit matrix}, a condensed version
of the incidence matrix of~$\D$. Dembowski~\cite{PD58, PD68} proved that
orbit matrices of $2$-designs satisfy a system of quadratic equations.
These equations were used by Janko and Tran Van Trung for the construction
of symmetric designs in~\cite{JTVT85} and in many subsequent papers.
The method of tactical decompositions was also used for the construction
of Steiner $2$-designs~\cite{JS87, VK02, CDDRS19, CDD25} and generalized
to $t$-designs~\cite{KNP11, KNP14, AN15, KW23}.

When constructing designs with a prescribed group of automorphisms,
orbit matrices can reduce the number of possibilities that need to
be considered. The greater the number of orbits, the stronger this
reduction becomes. For a point-transitive group~$G$, say a
group of order~$91$ acting on $S(2,6,91)$ designs, an orbit matrix
has just one row and does not contain any useful information on the
distribution of the points among the block-orbits.

A different approach, which is not directly affected by the number of
point-orbits, is the so-called Kramer-Mesner method; see~\cite{KM76}
and~\cite[Chapter~9.2]{KO06}. The existence of $t$-$(v,k,\lambda)$
designs with large~$t$ and small~$v$ was established by this method
\cite{ML84, KLM85, KR86, KR90, BKLW95, BKLW98, RL01}. The method was
also adapted to $q$-analogs of designs~\cite{BKW18b}; see~\cite{BKW18a}
for the relevant definitions. It was used in~\cite{BEOVW16} to construct
the first nontrivial examples of Steiner designs over~$\F_2$.
The parameters $S_2(2,3,13)$ are still the only concrete instance
where $q$-analogs of Steiner designs are known to exist.

The action of a permutation group~$G$ on~$V$ splits the set
${V\choose t}$ of all $t$-subsets into orbits $\T_1,\ldots,\T_m$, and
the set ${V\choose k}$ into orbits $\K_1,\ldots,\K_n$. Let $a_{ij}$ be
the number of $k$-subsets $K\in \K_j$ containing a given $t$-subset
$T\in \T_i$. Since the orbits form a tactical decomposition, this number
does not depend on the choice of $T\in \T_i$. Now any $t$-$(v,k,\lambda)$
design invariant under~$G$ corresponds to a $0$-$1$ solution~$x$ of the
system of linear equations
\begin{equation}\label{KMsys}
A\cdot x = \lambda \boldsymbol{1},
\end{equation}
where $A=[a_{ij}]$ is the so-called \emph{Kramer-Mesner matrix} and
$\boldsymbol{1}$ is the all-ones vector. In general, the entries
of~$A$ are non-negative integers. However, a Steiner design has
$\lambda=1$, so any column with entries greater than~$1$ cannot be
part of a solution. We can delete all such columns from the
Kramer-Mesner matrix, i.e.\ discard the corresponding orbits~$\K_j$.

Let us turn to the example of $S(2,6,91)$ designs with
an automorphism of order~$91$. The stabilizer of any $6$-subset
is trivial, hence all $6$-subset orbits are of size $91$.
The total number of these orbits is ${91\choose 6}/91 = 7\,324\,878$.
A computation in GAP~\cite{GAP} shows that most of them cover a
$2$-subset more than once and therefore cannot be block orbits of
a Steiner $2$-design. The remaining $1\,774\,964$ \emph{good
orbits} cover each $2$-subset at most once. Equivalently,
every pair of elements of a good orbit intersects in at most one
point or, for a general Steiner $t$-design, in at most $t-1$ points.
We  don't actually build all orbits and discard the
``bad'' ones. Instead, we build only the good orbits using an
algorithm developed in~\cite{KV16, RVK19, KVK21} and available
in the GAP package \emph{Prescribed automorphism groups}~\cite{PAG}.
The algorithm was developed for quasi-symmetric designs,
i.e.\ designs with only two block intersection sizes~$x$ and~$y$.
Steiner $2$-designs are quasi-symmetric with $x=0$ and $y=1$, so
the algorithm can be applied directly. We have also modified
the algorithm to build good orbits for general Steiner $t$-designs.

The most time-critical part of the computation is finding
solutions of the Kramer-Mesner system~\eqref{KMsys}.  Solving
systems of linear equations over $\{0,1\}$ is a well-known
NP-complete problem. The fifth author's program
\emph{solvediophant}~\cite{AW98, AW21} was instrumental for this
task in~\cite{BKLW95, BKLW98, RL01} and in many other papers.
The special case of Steiner $t$-designs offers the advantage that
all entries of the Kramer-Mesner matrix~$A$ are also in $\{0,1\}$,
provided that only good $k$-subset orbits are considered. In this
case, finding solutions of~\eqref{KMsys} is an instance of the
\emph{exact cover} problem. Although this problem is still
NP-complete~\cite{GJ79}, for practical purposes it allows for
very fast backtracking approaches: Knuth's algorithms
\textit{dancing links}~\cite{DEK00,DEK20} and \textit{dancing
cells}~\cite{DEK25}. Implementations of both algorithms can be
downloaded from \cite{DEK23}. The dancing links algorithm was a
crucial step for the construction of $S_2(2,3,13)$ designs
in~\cite{BEOVW16}.

For cyclic $S(2,6,91)$ designs, the reduced Kramer-Mesner
matrix~$A$ is of size $45\times 1\,774\,964$. The system~\eqref{KMsys}
has exactly $120$ solutions, which can be computed within a few
days of CPU time on a fast desktop computer, e.g.\ with Knuth's
dancing cells program \emph{ssxcc}. The result was independently
verified with Kaski and Pottonen's dancing links implementation
\emph{libexact}~\cite{KP08}. Many of the solutions correspond
to isomorphic designs. We used McKay and Piperno's programs
\emph{nauty} and \emph{Traces}~\cite{MP14} to sort out
non-isomorphic designs, and got exactly the four examples
of Table~\ref{tab1}. This agrees with~\cite[Table~1]{BT12}
and~\cite[Table~2]{BT19}, where it was reported that
there are exactly four cyclic difference families and
cyclically permutable perfect binary constant weight codes
for parameters $(91,6,1)$.


In Knuth's terminology, the rows of the $\{0,1\}$-matrix~$A$
are \emph{items} and the columns are \emph{options}. The exact
cover problem is to find sets of options covering every
item exactly once. The entries of~$A$ indicate which items are
covered by which options. For our example problem, the items are
$2$-subset orbits and the options are good $6$-subset orbits
under the action of the group of order~$91$. An $S(2,6,91)$
design has $b=\frac{v(v-1)}{k(k-1)}=273=3\cdot 91$ blocks, so
only three options need to be chosen. This is why we can solve
a problem with as many as $1\,774\,964$ options: the search
tree is very ``flat'' in this case. Usually, the depth of the
search tree grows proportionally with the number of options,
leading to an exponential increase of execution time.

There is another obstacle when applying this method to
large problems. We are interested in the number of
designs up to isomorphism. As we can see in our example,
many solutions of the system~\eqref{KMsys} correspond to
isomorphic designs. In a sense, we are performing equivalent
parts of the computation repeatedly. To prevent this, the
normalizer $N(G)=\{\pi \in S_v \mid \pi G = G\pi\}$ of the
prescribed group~$G$ in the symmetric group~$S_v$
can be used. If~$\D$ is a solution, i.e.\ a $G$-invariant
design ($\D^G=\D$), then for every $\pi \in N(G)$
the isomorphic design $\D^\pi$ is also a solution:
$$(\D^\pi)^G=\D^{\pi G}=\D^{G\pi}=(\D^G)^\pi = \D^\pi.$$
There are $n$ choices for the first orbit in an $m\times n$
Kramer-Mesner system~\eqref{KMsys}. If the orbits $\K_{j_1}$
and $\K_{j_2}$ are equivalent under the normalizer, i.e.\
there exists a permutation $\pi \in N(G)$ such that
$\K_{j_2}=\K_{j_1}^\pi$ holds, then these two choices will
lead to isomorphic designs. Therefore, we can restrict the
choice of the first orbit to a set~$\N$ of representatives
under the action of~$N(G)$. This trick has been used for
general $t$-$(v,k,\lambda)$ designs~\cite{PK05} and also
in the $q$-analog case~\cite{KKW18}.

It is known that the normalizer $N(C_v)$ of the cyclic
group~$C_v$ in~$S_v$ is $C_v\rtimes \Aut(C_v)$. Hence, for our
example problem, the order of $N(C_{91})$ is $91\cdot \phi(91) =
91\cdot 72 = 6552$, where $\phi$ is Euler's totient function.
Using the GAP package \emph{images}~\cite{JPWJ24}, we compute the
set~$\N$ of $N(G)$-representatives, reducing the number of
choices for the first orbit from $1\,774\,964$ to $|\N|=24\,717$.
The computation of $\mathcal{N}$ takes about an hour of CPU time,
but we shall soon see it is well worth the effort.

The straightforward way of applying the normalizer to~\eqref{KMsys}
would be to fix an orbit from~$\N$, modify the Kramer-Mesner
system so that it must be chosen (e.g.\ by adding an
extra equation), and run the solver. In the second run, another
orbit from~$\N$ is fixed and all orbits equivalent to the first
fixed orbit are discarded. Isomorphic designs containing one of
these orbits would already have been found in the first run.
In this way, the system is restricted further in every new run of
the solver. Our problem has a very large Kramer-Mesner matrix,
so the overhead for modifying the system and passing it to
the solver is considerable. Instead of restarting the
computation $24\,717$ times, we have found ways to perform
this strategy in a single run of the solver.

In a first approach, we add a copy of every column corresponding
to the $\N$-orbits to~$A$, as well as an extra row with $1$-entries
in place of the copies and on the right-hand side. The additional row
ensures that each solution must contain an orbit from~$\N$. The
reason for copying the $\N$-columns instead of putting $1$-entries
in the original $\N$-positions is not to miss solutions
involving several $\N$-orbits.

Going back to our example problem, we have now enlarged the
Kramer-Mesner matrix by one row and $24\,717$ columns,
thereby achieving that the solutions of the new system must contain
at least one $\N$-orbit. Solutions containing more than one
$\N$-orbit would be returned multiple times, with a $1$-entry
either in the original column or its copy. As it turned out,
there are no such solutions for this particular problem.
Using one of the above mentioned exact cover solvers, one
can find all solutions of the enlarged system in just a few
hours of CPU time. There are now only $8$ solutions: the
designs $\D_1$ and $\D_2$ are found once, and the designs
$\D_3$ and $\D_4$ are found three times each.

For the second approach, we focus on the observation that, after
having fixed an $\N$-orbit for a specific run, all orbits equivalent
to the fixed orbit can be discarded in subsequent runs. Indeed,
this can be modeled by a single instance of \emph{exact covering with
colors} as introduced by Knuth~\cite[pp.~85-89]{DEK20}. In this
extension to the exact cover problem, additional \textit{secondary
items} are allowed that have to be covered at most once. Moreover,
options may assign colors to the secondary items they cover.
In a solution, the secondary items may be covered by more than one
option, provided that the involved colors are identical.

In our case, we add a secondary item for each orbit $\K_i \in \N$
and assign this item to all options corresponding to the
$N(G)$-equivalent orbits, setting its color to~$1$.  Then, for
all copies of representatives $\K_j\in\N$ with $j>i$, we add the
same item with color~$0$. This ensures that, as soon as a
representative $\K_{i+1}\in \N$ is fixed, all orbits equivalent
to $\K_j\in \N$ with $j\leq i$ are discarded.  In
Section~\ref{sec3}, we give benchmark results for several
instances coming from $S(2,6,91)$ designs with prescribed groups
of order~$84$.

\section{Proof of Theorem~\ref{tm1}}\label{sec3}

In the design construction business, choosing a suitable
prescribed group of automorphisms is somewhat of an art.
In~\cite{CDD25}, a non-abelian group of order~$26$ was
prescribed, which appears as a group of automorphisms
of the McCalla design~$\D_2$ and the Colbourn design~$\D_3$.
The search revealed only these two already known designs.
We decided to avoid automorphisms of order~$13$ and prescribe
groups of order~$84=2^2\cdot 3\cdot 7$. The McCalla design
is invariant under a group $C_7\rtimes C_{12}$ of order~$84$.
It acts on the points in an orbit of size~$7$ and another
orbit of size~$84$. To be more systematic, we will consider
all groups of order~$84$ acting on the points of $S(2,6,91)$
designs in this way.

\begin{table}[!b]
\caption{The groups of order $84$ and results of the classification.}
\label{tab2}
\begin{tabular}{ccrrrc}
\hline
Group & Structure & \multicolumn{1}{c}{Orbits} & \multicolumn{1}{c}{$|N(G_i)|$} &
\multicolumn{1}{c}{$|\N|$} & Designs\\
\hline
$G_1$ & $C_7 \rtimes C_{12}$ & $703\,591$ & $7\,056$ & $8\,509$ & $8$ \\
$G_2$ & $C_4 \times (C_7 \rtimes C_3)$ & $637\,595$ & $7\,056$ & $7\,697$ & $8$ \\
$G_3$ & $C_7 \times (C_3 \rtimes C_4)$ & $757\,275$ & $42\,336$ & $8\,985$ & $0$ \\
$G_4$ & $C_3 \times (C_7 \rtimes C_4)$ & $883\,955$ & $14\,112$ & $5\,443$ & $0$ \\
$G_5$ & $C_{21} \rtimes C_4$ & $1\,279\,623$ & $42\,336$ & $2\,697$ & $0$ \\
$G_6$ & $C_{84}$ & $1\,011\,339$ & $14\,112$ & $35\,765$ & $0$ \\
$G_7$ & $C_2 \times (C_7 \rtimes C_6)$ & $30\,191$ & $7\,056$ & $406$ & $0$ \\
$G_8$ & $S_3 \times D_{14}$ & $2\,443$ & $21\,168$ & $23$ & $0$ \\
$G_9$ & $C_2 \times C_2 \times (C_7 \rtimes C_3)$ & $378\,903$ & $21\,168$ & $1\,593$ & $2$ \\
$G_{10}$ & $C_7 \times A_4$ & $409\,764$ & $84\,672$ & $2\,018$ & $0$ \\
$G_{11}$ & $(C_{14} \times C_2) \rtimes C_3$ & $577\,269$ & $42\,336$ & $1\,184$ & $6$ \\
$G_{12}$ & $C_6 \times D_{14}$ & $61\,021$ & $14\,112$ & $444$ & $0$ \\
$G_{13}$ & $C_{14} \times S_3$ & $278\,489$ & $42\,336$ & $2\,184$ & $0$ \\
$G_{14}$ & $D_{84}$ & $4\,265$ & $42\,336$ & $94$ & $0$ \\
$G_{15}$ & $C_{42} \times C_2$ & $666\,585$ & $42\,336$ & $7\,162$ & $0$ \\
\end{tabular}
\end{table}

According to the GAP library of small groups~\cite{GAP},
there are $15$ non-isomorphic groups of order~$84$.
We list them in Table~\ref{tab2} with a description of
the structure, numbered as they appear in the GAP library.
Using GAP, it is easy to check that each of the $15$ groups
has a single subgroup of order~$12$ up to conjugation.
This means that it allows a unique transitive action
on~$7$ points, and this is also true for the orbit of
size~$84$. We conclude that each of the groups has
a single permutation representation on two orbits of
sizes~$7$ and~$84$ up to permutation isomorphism, i.e.\
up to conjugation in the symmetric group~$S_{91}$.

For each of these permutation groups, we perform a computer
classification of $S(2,6,91)$ designs as described in the
previous section. Results of the calculations are
reported in Table~\ref{tab2}. The column labeled `Orbits'
contains the number of good orbits of $6$-subsets of~$V$,
calculated in PAG~\cite{PAG}. The next column contains the
order of the normalizer, calculated by the GAP~\cite{GAP}
command \texttt{Normalizer}. The column labeled $|\N|$
contains the number of $\N$-representatives of $6$-subset
orbits, calculated using the \emph{images}~\cite{JPWJ24}
package. Finally, the last column contains the number
of non-isomorphic designs obtained from each group.

We used two exact cover programs to solve
each Kramer-Mesner system, namely \emph{ssxcc}~\cite{DEK25}
and \emph{libexact}~\cite{KP08}, and the results always
agreed. Over the $15$ groups of order~$84$, the required
CPU time varied from a few hours to just a few seconds. The
latter cases include~$G_8$ and~$G_{14}$, with few orbits
and no solutions. Interestingly, the groups~$G_5$ and~$G_6$
with most orbits required only a few minutes of CPU time.
Longer running times were required for the cases admitting
solutions. The numbers of non-isomorphic designs and their
full automorphism groups were calculated with \emph{nauty}
and \emph{Traces}~\cite{MP14}.

To conclude, four of the $15$ groups of order~$84$ give rise
to designs: $G_1$, $G_2$, $G_9$ and $G_{11}$. As expected,
the McCalla design~$\D_2$ was recovered from~$G_1$.
This is the only case where the full automorphism group
is larger than the prescribed group. All other constructed
designs have full automorphism groups of order~$84$.
This also means that designs coming from different groups
cannot be isomorphic. The sum of the numbers in the last
column of Table~\ref{tab2} is $24$, proving Theorem~\ref{tm1}.
The constructed designs are available on our web page in
GAP-compatible format:
\begin{center}
\url{https://web.math.hr/~krcko/results/steiner.html}
\end{center}

\begin{table}[t]
    \caption{Benchmark results.}
    \label{tab3}
    \begin{tabular}{cclrrrr}
        \hline
        Group    & Method & Items   & Options  & Sols. & Nodes   & Time (s)  \\
        \hline
        $G_1$    & a)    & $50+0$    & $703\,591$ & $672$       & $324\,581$  & $10\,914$ \\
                 & b)    & $51+0$    & $712\,100$ & $56$        & $156\,630$  & $5\,705$  \\
                 & c)    & $51+8508$ & $712\,100$ & $8$         & $72\,694$   & $5\,016$  \\
        \hline
        $G_2$    & a)    & $50+0$    & $637\,595$ & $672$       & $388\,315$  & $8\,176$  \\
                 & b)    & $51+0$    & $645\,292$ & $56$        & $161\,629$  & $3\,338$  \\
                 & c)    & $51+7696$ & $645\,292$ & $8$         & $60\,911$   & $2\,796$  \\
        \hline
        $G_9$    & a)    & $51+0$    & $378\,903$ & $504$       & $702\,346$  & $8\,633$  \\
                 & b)    & $52+0$    & $380\,496$ & $43$        & $115\,548$  & $1\,267$  \\
                 & c)    & $52+1592$ & $380\,496$ & $2$         & $19\,933$   & $792$   \\
        \hline
        $G_{11}$ & a)    & $51+0$    & $577\,269$ & $3\,024$    & $2\,390\,706$ & $19\,322$ \\
                 & b)    & $52+0$    & $578\,453$ & $241$       & $255\,622$  & $2\,611$  \\
                 & c)    & $52+1183$ & $578\,453$ & $6$         & $32\,015$   & $1\,084$  \\
    \end{tabular}
\end{table}

The run times for solving the Kramer-Mesner systems for $G_1$, $G_2$,
$G_9$, and $G_{11}$ are compared in Table~\ref{tab3}. We solved each problem
a)~without using the normalizer,
b)~using the normalizer, but without discarding orbits, and
c)~using the normalizer and discarding orbits. For each instance, we list the number of
items (primary$+$secondary), the number of options, the number of solutions,
the number of nodes of the search tree, and the computing times in seconds
using \textit{ssxcc}. The number of secondary items is one less than~$|\N|$.
It can be seen that the search tree becomes much smaller if the normalizer is
used, and it becomes even smaller if orbits are discarded. However, the latter approach
comes with a huge overhead if there are many orbits, so that the actual run time
may not benefit much. If the number of inequivalent orbits under the action of
the normalizer is reasonable, i.e.~less than $2\,000$, the discarding of orbits
outperforms the other approaches.

\section{Concluding remarks}\label{sec4}

In this paper, we have fine-tuned the Kramer-Mesner approach
for constructing Steiner $2$-designs with prescribed
groups of automorphisms. We have used it to construct new
$S(2,6,91)$ designs with groups of orders~$84$. Software
to perform the calculations is already partially available in
the GAP package PAG~\cite{PAG}, and we are planning to include
the new improvements in future releases.

An obvious research direction is to try constructing designs
with admissible parameters $S(2,k,v)$ for which the existence
question is open, such as the ones mentioned in the Introduction.
The challenge lies in finding suitable permutation groups that
can act as groups of automorphisms of the designs and are large
enough so that the necessary calculations can be performed.
Presumably, many groups have already been tried for the
parameters that are still open. It might be difficult to
find just the right groups to get designs, if they exist.

A more accessible goal could be to construct further examples
with parameters for which a few designs are already
known, as we have done for $S(2,6,91)$. The parameters
$S(2,7,91)$ are similar, and there are only two known
examples, constructed in~\cite{BNKKS83}. The full automorphism
groups of these two designs are both of order $1092$ and they
are not isomorphic. For each of the groups in Table~\ref{tab2},
we have also performed a classification of $S(2,7,91)$ designs,
but no new examples were found. The two known designs were
recovered from groups~$G_1$ and~$G_2$.

A better choice of prescribed automorphism group might lead
to new examples of $S(2,7,91)$ designs. Other parameters for
which only a few designs are known include $S(2,6,v)$ for $v=66$,
$76$, $96$, $106$, and $111$. The first examples of these
designs were constructed by Denniston~\cite{RHFD80} and
Mills~\cite{WHM75b, WHM78, WHM79}.

\vskip 2mm

\textsc{Acknowledgements.}

This work was done while the second author was visiting
the University of Bayreuth supported by a scholarship from the German
Academic Exchange Service (DAAD). The second and the fourth authors
have been partially supported by the Croatian Science Foundation
under the project $9752$.

\end{document}